\documentclass[a4paper, draft]{article}
\usepackage{amsmath, amssymb, amsthm, url}
\title{Spectral expansions of overconvergent modular functions}
\author{David Loeffler \\
Department of Mathematics \\ 
Imperial College \\ 
South Kensington, London SW7 2AZ, UK\\
{\tt david.loeffler@imperial.ac.uk}}
\newtheorem{theorem}{Theorem}
\newtheorem{corollary}[theorem]{Corollary}
\newtheorem{lemma}[theorem]{Lemma}
\newtheorem{conjecture}[theorem]{Conjecture}
\DeclareMathOperator{\Res}{Res}
\DeclareMathOperator{\Ker}{Ker}
\DeclareMathOperator{\SL}{SL}
\begin{document}
\maketitle

\begin{abstract}
The main result of this paper is an instance of the conjecture made by Gouv\^ea and Mazur in \cite{GMsearch}, which asserts that for certain values of $r$ the space of $r$-overconvergent $p$-adic modular forms of tame level $N$ and weight $k$ should be spanned by the finite slope Hecke eigenforms. For $N = 1$, $p = 2$ and $k = 0$ we show that this follows from the combinatorial approach initiated by Emerton \cite{emerton} and Smithline \cite{smithline-thesis}, using the classical $LU$ decomposition and results of Buzzard--Calegari \cite{BCslopes}; this implies the conjecture for all $r \in (\frac{5}{12}, \frac{7}{12})$. Similar results follow for $p = 3$ and $p = 5$ with the assumption of a plausible conjecture, which would also imply formulae for the slopes analogous to those of \cite{BCslopes}.

We also show that (for general $p$ and $N$) the space of weight 0 overconvergent forms carries a natural inner product with respect to which the Hecke action is self-adjoint. When $N = 1$ and $p \in \{2, 3, 5, 7, 13\}$, combining this with the combinatorial methods allows easy computations of the $q$-expansions of small slope overconvergent eigenfunctions; as an application we calculate the $q$-expansions of the first 20 eigenfunctions for $p = 5$, extending the data given in \cite{GMsearch}.
\end{abstract}
\section{Background}

Let $S_k(\Gamma_1(N))$ denote the space of classical modular cusp forms of weight $k$ and level $N$. It has long been known that these objects satisfy many interesting congruence relations. One very powerful method for studying the congruences obeyed by modular forms modulo powers of a fixed prime $p$ is to embed this space into the $p$-adic Banach space $\mathcal{S}_k(\Gamma_1(N),r)$ of $r$-overconvergent $p$-adic cusp forms, defined as in \cite{katz} using sections of $\omega^{\otimes k}$ on certain affinoid subdomains of $X_1(N)$ obtained by removing discs of radius $p^{-r}$ around the supersingular points; this space has been used to great effect by Coleman and others (\cite{C-CO, Cpadic}).

It is known that there is a Hecke action on $\mathcal{S}_k(\Gamma_1(N),r)$, as with the classical spaces, and these operators are continuous; and moreover, at least for $0 < r < \frac{p}{p+1}$, the Atkin-Lehner operator $U$ is compact. There is a rich spectral theory for compact operators on $p$-adic Banach spaces (see \cite{serre-endo}), and this is a powerful tool for studying the spaces $\mathcal{S}_k(\Gamma_1(N), r)$. In this paper, we shall attempt to make this spectral theory explicit in the case $N = 1$, $k = 0$, for certain small primes $p$.

\section{A useful basis}

In all the computations in this paper, we shall restrict to the case of tame level 1; hence we shall write $\mathcal{S}_k(r)$ for $\mathcal{S}_k(\SL_2(\mathbb{Z}), r)$, regarded as a Banach space over $\mathbb{C}_p$.

Recall that if $\psi$ is any lifting of the mod $p$ Hasse invariant to a modular form in characteristic $0$, and $E$ is any elliptic curve over $\mathbb{C}_p$ such that $|\psi(E)| > p^{-\frac{p}{p+1}}$, then $E$ has a canonical $p$-subgroup; hence, for $0 < r < \frac{p}{p+1}$, the $r$-overconvergent locus $X_0(1)_{\ge p^{-r}}$ is isomorphic to a certain subregion of $X_0(p)$. (This is proved in \cite{katz}, using the theory of the Newton polygon.)

If $p$ is one of the primes 2, 3, 5, 7, or 13, then $X_0(p)$ has genus 0. We shall pick an explicit uniformiser for this curve, and identify in terms of this uniformiser the image of $X_0(1)_{\ge p^{-r}}$ under the canonical subgroup map, and hence obtain a basis for our space $\mathcal{S}_k(r)$. 

\begin{theorem}\label{basis}
Let $p$ be one of the primes 2, 3, 5, 7, or 13. Let $f_p$ be the function
\[\left[\frac{\Delta(pz)}{\Delta(z)}\right]^{\frac{1}{p-1}}.\]
Then $f_p$ is a rational function on the modular curve $X_0(p)$, and the forgetful functor gives an isomorphism between the region of the modular curve $X_0(p)$ where $|f_p| \le 1$ and the ordinary locus $X_0(1)_{\rm ord}$. Moreover, for any $r \in [0, \frac{p}{p+1})$, this extends to an isomorphism between the region where $|f_p| \le p^{\frac{12r}{p-1}}$ and $X_0(1)_{\ge p^{-r}}$. 
\end{theorem}

\begin{proof}
That $f_p$ is a rational function on $X_0(p)$ is clear from the fact that $\Delta(z)$ and $\Delta(pz)$ are both classical modular forms of weight 12 and level $p$, and $\Delta$ has no zeros on $X_0(p)$. It has a zero of order 1 at $z=\infty$ by inspection of its $q$-expansion, and no other zeros as $\Delta$ does not vanish on the complex upper half-plane; so it is a uniformiser for $X_0(p)$.

It remains to prove that the subsets defined by $|f_p| \le p^{\frac{12r}{p-1}}$ agree with the $r$-overconvergent locus as defined in \cite{katz} using lifts of the Hasse invariant. For $p = 2$ this is proved in \cite[\S 4]{BCslopes}; for $p \ge 5$ it is \cite[Prop 3.5]{smithline-bounding}. In the remaining case $p=3$ Smithline uses a different measure of supersingularity and it is not immediately obvious this agrees with the valuation of the Hasse invariant; we show that the two do in fact agree below, in \S \ref{appendix}.
\end{proof}

\begin{corollary}
For any $0 \le r < \frac{p}{p+1}$, the space $\mathcal{S}_0(r) = \mathcal{S}_0(\SL_2(\mathbb{Z}), r)$ of $r$-overconvergent $p$-adic tame level 1 cuspidal modular functions (modular forms of weight 0) has an orthonormal basis $(cf_p, (cf_p)^2, (cf_p)^3, \dots)$ where $c$ is any element of $\mathbb{C}_p$ with $|c| = p^{\frac{12r}{p-1}}$.
\end{corollary}

(This follows as we have given an isomorphism between this space and a $p$-adic closed disc, and the algebra of rigid-analytic functions on a $p$-adic closed disc with uniformising parameter $x$ is the Tate algebra $\mathbb{C}_p\langle x \rangle$.)

\begin{theorem}\label{recurrence}
Let $U$ be the Atkin-Lehner operator acting on $\mathcal{S}_0(r)$, and let $u_{ij}^{(r)}$ be the matrix coefficients of $U$ with respect to the basis defined above. Then the following results hold:
\begin{enumerate}
\item $u_{ij}^{(r)} = c^{j-i} u_{ij}^{(0)}$.
\item There is a $p \times p$ matrix $M^{(r)}$, which is `skew upper triangular' (that is, $M_{ij}^{(r)} = 0$ if $i + j > p + 1$), with the property that 
\[u_{ij} = \sum_{a,b = 1}^p M^{(r)}_{ab} u_{i-a, j-b}^{(r)}\]
for all $i,j > p$.
\item $u_{ij}^{(r)} = 0$ if $i > pj$ or $j > pi$, so in particular $U(f_p^k)$ is a polynomial in $f_p$ of degree at most $pk$.
\end{enumerate}
\end{theorem}

\begin{proof}
Part (1) is an elementary manipulation. Given this, it is clearly sufficient to prove the existence of $M$ when $r = 0$. This result is well-known for $p=2$, and may be found in Emerton's thesis \cite{emerton}; it is apparently initially due to Kolberg. The same approach may be used for the other values of $p$, or alternatively one may deduce the result from \cite[Lemma 3.3.2]{smithline-thesis}, where it is shown that there is a polynomial $I_p(x,y)$ of degree $p$ in each variable such that $I_p(V(f_p), \frac{1}{f_p})=0$, where $V$ is the operator induced by $q \mapsto q^p$. Smithline produces this identity by noting that there exists a polynomial $H_p$ of degree $p+1$ with integer coefficients such that $\frac{H_p(f_p)}{f_p}$ is the level 1 $j$-invariant, and thus we have
\[\frac{H_p(p^{-12/(p-1)}/f_p)}{p^{-12/(p-1)}/f_p} = \frac{H_p(V(f_p))}{V(f_p)}\]
since both sides are equal to $V(j)$. Clearing denominators and cancelling the factor $V(f_p) - p^{-12/(p-1)}/f_p$ (which is clearly not identically zero) gives $I_p$, and it is thus clear that $I_p$ has integer coefficients, total degree $p + 1$, constant coefficient equal to 1 and all linear terms zero. Multiplying by $f_p^j$, applying $U$ and using ``Coleman's trick'' --- the identity $U(fV(g)) = gU(f)$ --- gives the required recurrence, with $M_{ab}$ being the coefficient of $x^a y^b$ in $-I_p(x,y)$. So part (2) of the theorem follows. 

Finally, since $U(1) = 1$ and coefficients of the recurrence are polynomials in $f_p$ of degree at most $p$, it follows by induction that $U(f_p^j)$ must be a polynomial of degree at most $pj$ in $f_p$; thus $u_{ij} = 0$ if $i > pj$. On the other hand, it is immediate from the $q$-expansion that if $j > pi$, $U(f_p^j)$ must vanish to degree $i$ at the origin, so $u_{ij} = 0$ in this region as well.
\end{proof}

The polynomials $H_p$ are easy to compute by comparing $q$-expansions, and hence we can easily determine the polynomials $I_p$ explicitly (they are tabulated in \cite[\S 3.3]{smithline-thesis}) and thus the matrices $M$. For example, when $p=2$ we find that
\[M^{(0)} = \begin{pmatrix} 48 & 1 \\ 2^{12} & 0 \end{pmatrix},\]
and when $p = 3$,
\[M^{(0)} = \begin{pmatrix} 270 & 36 & 1\\
26244 & 729 & 0\\
531441 & 0 & 0
\end{pmatrix}.\]

\begin{corollary}\label{corollary}
The operator $U$ is an ``operator of rational generation'' in Smithline's sense; that is, there exists a rational function $R(x,y)$ whose Taylor series expansion is equal to $\sum_{i,j} u_{ij} X^i Y^j$. The function $R$ is equal to
\[-\frac{y}{p} \frac{\partial}{\partial y} \log I_p(x,y).\]
\end{corollary}

\section{Computations of slopes}

If $X$ is any compact operator acting on a $p$-adic Banach space, it has a (possibly empty!) countable set of nonzero eigenvalues, for each of which the generalised eigenspace $\bigcup_{k = 1}^\infty \Ker\left[ (U-\lambda_i)^k\right]$ is finite-dimensional. The $p$-adic valuations of these eigenvalues are known as the {\it slopes}. The finite slope eigenvalues occur as the inverses of roots of the characteristic power series $\det (I - tX)$. 

In our case, it is known that $U$ is compact for $r \in (0, \frac{p}{p+1})$. Given the values of $u_{ij}^{(r)}$ for $1 \le i,j \le N$, it is easy to calculate the characteristic power series of this $N \times N$ matrix (since the entries are rational); and the general theory of compact operators tells us that this will converge rapidly to the characteristic power series of $U$. So we can easily calculate approximations to the eigenvalues, and in particular we can determine the slopes. The results obtained will be independent of $r$, since it is known that any overconvergent $U$-eigenform of finite slope must extend to a function on $X_0(1)_{\ge p^{-r}}$ for all $r < \frac{p}{p+1}$ (see \cite{Banal}).

The slopes of $U$ are somewhat mysterious; the complete list of slopes is known only for $p=2$, tame level 1 and weight 0 by \cite{BCslopes}, and for 2-adic, 3-adic and 5-adic weights near the boundary of weight space by \cite{BKbound}, \cite{jacobs} and \cite{kilford5} respectively. There are conjectures (\cite{Bques}, \cite{clay}) for a general weight, prime and level, but these appear to be rather inaccessible at present.

In the approach of \cite{BCslopes}, the next step would be to attempt to decompose the $U$ operator as $U = ADB$ where $A$ is lower triangular, $B$ is upper triangular, $D$ is diagonal, and both $A$ and $B$ have all diagonal entries 1. If this factorisation exists (which is the case if none of the top left $r \times r$ minors are singular) then it is unique, and can be calculated rapidly by Gaussian elimination; usefully, the $i, j$ entry of each of $A,B,D$ is determined by $u_{mn}$ for $mn \le \max(i,j)$, so in our case the entries of these matrices are rational and can be calculated exactly using our algorithm for calculating $U$.

\begin{conjecture}\label{conjecture}
For $p \in \{2, 3, 5\}$ and all $r$ in some open interval containing $\frac{1}{2}$, the $U$ operator acting on $\mathcal{S}_k(r)$ has a factorisation $U^{(r)} = A^{(r)}DB^{(r)}$, where $A^{(r)}$ and $B^{(r)}$ have entries in $\mathcal{O}_{\mathbb{C}_p}$ and are congruent to the identity modulo $p$, and the entries of $D$ are given by the following formulae:
\[
\begin{array}{r|c|c|}
p & D_{ii} & \nu_p D_{ii}\\
\hline
2 & \dfrac{2^{4i + 1} (3i)!^2 i!^2}{3\cdot (2i)!^4} & 1 + 2 \nu_2\left(\frac{(3i)!}{i!}\right)\\
3 & \dfrac{3^{3i} (6i)!(2i)!i!}{2\cdot(3i)!^3} & 2i + 2 \nu_3 \left( \frac{(2i)!}{i!}\right)\\
5 & \dfrac{5^{2i} (10i)!(3i)!^2i!}{3\cdot(5i)!^3(2i)!} & i + 2\nu_5\left(\frac{(3i)!}{i!}\right)\\
\end{array}
\]
\end{conjecture}

This is known in the case $p = 2$, by \cite{BCslopes} (for $r = \frac{1}{2}$, but we extend the result to all $r \in (\frac{5}{12}, \frac{7}{12})$ below). For $p = 3$ and $p = 5$ it is open, but a calculation of $U_{ij}$ for $1 \le i, j \le 100$ suggests that the conjecture holds for $r \in (\frac{1}{3}, \frac{2}{3})$ in both cases. However, the same computation suggests that the entries of $A$ and $B$ are not given by any hypergeometric term (as they are divisible by too many large primes).

If this conjecture is true, then lemma 5 of \cite{BCslopes} would tell us that the Newton polygon of $ADB$ is the same as that of $D$, so the $i$th slope would be equal to the valuation of the $i$th diagonal entry of $D$. Indeed, Frank Calegari has conjectured formulae for the slopes for $p=3$ and $p=5$ (cited in \cite{smithline-ratgen}), and these agree with those given in the third column above. Furthermore, these formulae also appear to agree with the combinatorial recipe of \cite{Bques}; but without a concise formula for $A_{ij}$ and $B_{ij}$, there does not seem to be any chance of proving these results by this method.

For $p = 7$ and $p = 13$ the pattern is much less clear; there still appears to be an $ADB$ factorisation with $A$ and $B$ congruent to the identity, but the entries of $D$ do not appear to be given by any simple hypergeometric form. It is interesting to note that in these cases, there are several distinct ``slope modules'' in the conjectural picture of \cite{clay}, so one would not expect all the slopes to be given by a single simple formula.

\section{Computations of eigenfunctions}

If $M$ is an $n \times n$ matrix over a $p$-adic field, then calculating the eigenvalues and eigenvectors of $M$ to any desired degree of accuracy is computationally very easy, as Hensel's lemma allows easy calculation of the eigenvalues. More generally, if $M$ is the matrix of a compact operator and $M_n$ is the $n \times n$ truncation, then one can calculate the eigenvectors of $M$ using $M_n$: if $\lambda$ is an eigenvalue of $M$, and $n$ is sufficiently large compared to the slope of $\lambda$, then there will be an eigenvalue $\lambda_n$ of $M_n$ which is highly congruent to $\lambda$, and and as $n \to \infty$, $\lambda_n$ will converge to $\lambda$ and the associated eigenvectors $v_n$ will converge to an eigenvector of $M$.

Let us do this in the case $p = 5$ (for comparison with the calculations in \cite{GMsearch}). We begin by fixing a value of $r$; in this case, it is convenient to choose $r = \frac{1}{3}$, since in this case we may take $c = p$ and the $u_{ij}$ are all rational. We now take an $N \times N$ truncation of the matrix of $U$ and diagonalise this using the PARI/GP functions \texttt{polrootspadic()} and \texttt{matker()}; this gives an approximate $U$-eigenfunction. As it is necessary to divide by entries of the matrix in this computation, the resulting eigenvector is known to slightly less precision than the eigenvalue; but this is not a serious problem as calculating the roots of $p$-adic polynomials is computationally very easy -- working modulo $5^{300}$ is no problem on current machines.

If we take $N = 3$, we obtain three eigenvalues of slopes $\sigma_1 = 1$, $\sigma_2 = 4$ and $\sigma_3 = 5$, and three corresponding approximate eigenfunctions $\phi_1$, $\phi_2$ and $\phi_3$. Repeating the calculation for a range of $N$, it seems that changing $N$ does not change $\phi_1 \bmod 5^{8}$, so the value obtained for $N=3$ is apparently already correct to this precision; moreover, taking $N=4$ is enough to give it mod $5^{10}$, and $N = 5$ gives it mod $5^{16}$. So the functions obtained appear to be converging very rapidly in the $q$-expansion topology (or, equivalently, in the supremum norm on $X_0(1)_{\rm ord}$). The first 30 terms of the $q$-expansion of the first few $\phi_i$ is given modulo $5^{15}$ in \S \ref{tables}.

\section{Spectral expansions}
\label{spectralexpansions}

It is a standard consequence of the spectral theory that for each nonzero eigenvalue $\lambda_i$ of $U$, there is a projection $\pi_i$ onto the corresponding generalised eigenspace, and this projection commutes with $U$. Since for any $x \ge 0$, the set $\Lambda_x$ of eigenvalues of slope $\le x$ is finite, one can form for any $h \in S_k(r)$ the series
\[ e_x(h) = \sum_{\lambda \in \Lambda_x} \pi_\lambda(h).\]

This is known as the {\it asymptotic $U$-spectral expansion of $h$}. This will not generally converge as $x \to \infty$; but it is uniquely determined by the property that for any $x$ there exists $\epsilon > 0$ with $\nu_p\left(\left\|X^k(h - e_x(h) )\right\|\right) \ge (x + \epsilon)k$ for all $k \gg 0$.

For $p = 2, 3, 5$, all the eigenspaces are conjecturally one-dimensional, spanned by eigenfunctions $\phi_i$, so we should obtain a sequence of constants $c_i(h) = \pi_i(f) / \phi_i$. In principle, the spectral theory gives an explicit form for the spectral projections $\pi_i$. The first projection $\pi_1$ is easy, as one simply iterates the process of applying $U$ and dividing by the eigenvalue $\lambda_1$. One can then consider $h' = h - \pi_1(h)$ and iterate $U$ on this; the same process of iterating and dividing by $\lambda_2$ should converge to the second projection $\pi_2$, but this is unstable with regard to small errors in the calculation of $\pi_1(h)$ -- such errors will inevitably grow at a rate of $(\lambda_1 / \lambda_2)^k$ until they swamp the desired answer. So this method is not really usable in practice.

However, the symmetry properties of $U$ provide us with an alternative approach. Let $g = p^{6/(p-1)}f$, so $(g, g^2, g^3, \dots)$ are a basis for $\mathcal{S}_0(\frac{1}{2})$.

\begin{theorem}
Define the symmetric bilinear form $\langle,\rangle$ on $\mathcal{S}_0(\frac{1}{2})$ by 
\[ \langle g^i, g^j \rangle = 
\begin{cases}
i & (i=j)\\
0 & (i \ne j)
\end{cases}.\]
Then $U$ is self-adjoint with respect to this form; and for all $i$ such that the $\lambda_i$ eigenspace is 1-dimensional and $\langle \phi_i, \phi_i \rangle \ne 0$, the spectral projection operators $\pi_i$ are given by 
\[\pi_i(h) = \frac{\langle h, \phi_i\rangle}{\langle \phi_i, \phi_i\rangle}.\]

Furthermore, the same formula in fact gives us a pairing $\mathcal{S}_0(r) \times \mathcal{S}_0(1-r) \to \mathbb{C}_p$ for any $r \in (\frac{1}{p+1}, \frac{p}{p+1})$.
\end{theorem}

\begin{proof}
If $p \in \{2, 3, 5, 7, 13\}$, then we can show that $U$ is self-adjoint with respect to this bilinear form by proving that $u_{ij}^{(1/2)} = \frac{j}{i} u_{ji}^{(1/2)}$. This follows from Corollary \ref{corollary} above; the generating function $R(x,y)$ is $\frac{y}{p} \frac{\partial}{\partial y}\log I_p(x,y)$, and from the construction of $I_p$ we see that it satisfies 
\[I_p(x,y) = I_p(p^{-12/(p-1)}y, p^{12/(p-1)}x),\]
so after an appropriate rescaling we see that $x \frac{\partial}{\partial x}R(x,y)$ is symmetric in $x$ and $y$, implying the result.

However, one can prove this in general -- without the assumption that $X_0(p)$ have genus 0 -- by using the theory of residues of $p$-adic differential forms. This theory is developed in \cite{FvdP}; for a general rigid space $X/k$ we can construct sheaves of finite differentials $\Omega_{X/k}^f$, and the notion of residue of a differential at a point can be defined in a consistent way. Now, if $\alpha$ and $\beta$ are in $\mathcal{S}_0(\frac{1}{2})$, and $w$ denotes the Atkin-Lehner involution on $X_0(p)$, then the differential
\[w^*(\alpha).\mathrm{d}\beta\]
is defined on the annulus $|A| = p^{-1/2}$ (a ``ring domain'') and thus has a residue at the cusp $\infty$. It is readily seen that if we define
\[\langle \alpha, \beta \rangle = \mathrm{Res}_{z = \infty} w^*(\alpha).\mathrm{d}\beta \]
then this agrees with the above definition when $p \in \{2, 3, 5, 7, 13\}$ (it is sufficient to check the result when $\alpha$ and $\beta$ are powers of $f$; in this case it is immediate from the fact that $w^*(g) = \frac{1}{g}$.)

Let $\Phi_1$ and $\Phi_2$ be the two canonical maps $X_0(p^2) \to X_0(p)$, namely $\Pi_1: (E, C) \mapsto (E, C[p])$ and $\Pi_2: (E, C) \mapsto (E/C[p], C/C[p])$; this gives a symmetric correspondence on $X_0(p)$, and the operator on functions corresponding to the trace of this correspondence is $U$. So we may write
\begin{align*}
\langle U\alpha, \beta \rangle &= \Res_{\infty \in X_0(p)} w^*(U\alpha)\ \mathrm{d} \beta\\
&= \Res_{\infty \in X_0(p)} U(w^* \alpha)\ \mathrm{d} \beta\\
&= \Res_{\infty \in X_0(p)} \Phi_{2*} \Phi_1^* w^*\alpha\ \mathrm{d} \beta \\
&= p \Res_{\infty \in X_0(p^2)} \Phi_{1}^* w^* \alpha\ \mathrm{d} \Phi_2^* \beta\\
&= \Res_{\infty \in X_0(p)} w^* \alpha\ d \Phi_{1*} \Phi_2^* \beta\\
&= \langle \alpha, U\beta \rangle.
\end{align*}

It now follows that any two eigenfunctions with different eigenvalues must be orthogonal, and the explicit form for the spectral projection operators is immediate.
\end{proof}

(Exactly the same argument also shows that the operators $T_\ell$ are self-adjoint for $\ell \ne p$.)

This pairing allows us to calculate spectral expansions extremely easily for functions $h$ that are at least $\frac{1}{2}$-overconvergent, given sufficiently accurate knowledge of the eigenfunctions themselves. As in the previous section, we shall take $p = 5$. Then the function $h = \frac{1}{j}$ is $r$-overconvergent for all $r < \frac{5}{6}$, and the constants $c_j$ turn out to be:
\[\begin{array}{|r|c|}
\hline
j & c_j \\
\hline
1 & 8295001 \\
2 & 5^{4}\times7540786 \\
3 & 5^{4}\times2165317 \\
4 & 5^{8}\times8075994 \\
5 & 5^{9}\times4502966 \\
6 & 5^{10}\times4930721 \\
7 & 5^{12}\times7120582 \\
8 & 5^{14}\times7314891 \\
9 & 5^{18}\times2324226 \\
10 & 5^{22}\times1076376 \\
\vdots & \vdots\\
\hline
\end{array}
\]
Here, as in the tables of eigenfunctions in \S \ref{tables}, we use a relative precision of $O(5^{10})$ -- that is, we write a general element of $\mathbb{Z}_5$ in the form $5^a b$ where $b \in (\mathbb{Z}/5^{10}\mathbb{Z})^\times$. These numbers appear to be tending 5-adically to zero extremely rapidly, suggesting that the $U$-spectral expansion is in fact convergent, at least in the (rather feeble) $q$-expansion topology.

One might optimistically make the following conjecture:

\begin{conjecture}[Gouv\^ea-Mazur spectral expansion conjecture, strong form] Let $h$ be any $r$-overconvergent modular function, where $r \in (\frac{1}{p+1}, \frac{p}{p+1})$. Then the spectral expansion of $h$ converges to $h$, in the supremum norm of $X_0(1)_{\ge p^{-r}}$.
\end{conjecture}

One cannot expect this to work for $r \le \frac{1}{p+1}$, for two reasons. Firstly, since the eigenfunctions themselves are not necessarily any more than $\frac{p}{p+1}$-overconvergent, we cannot guarantee that the linear functional $\langle \cdot, \phi_i\rangle$ even makes sense. More seriously, if $r < \frac{1}{p+1}$ then there exist nonzero functions in the kernel of $U$; the spectral expansion of any such form is always zero.

\section{The spectral expansion conjecture}

Let us now suppose either that $p = 2$, or that $p = 3$ or $5$ and Conjecture \ref{conjecture} above holds. We shall show that this implies the spectral expansion conjecture. 

Let $A^{(r)}$ and $B^{(r)}$ be the matrices occurring in the $LDU$ factorisation of $U^{(r)}$. ($D$ is clearly independent of $r$.)

\begin{lemma} For $p = 2$, Conjecture \ref{conjecture} holds for all $r \in \left(\frac{5}{12}, \frac{7}{12}\right)$; that is, for any $r$ in this range, $A^{(r)}$ and $B^{(r)}$ have entries in $\mathcal{O}_{\mathbb{C}_2}$ and their reductions modulo the maximal ideal are equal to the identity matrix.
\end{lemma}

\begin{proof}
Since by construction $A$ is lower triangular, $B$ is upper triangular and their diagonal entries are 1, it is sufficient to prove that $A^{(\frac{7}{12})}$ and $B^{(\frac{5}{12})}$ have entries in $\mathcal{O}_{\mathbb{C}_2}$. Conveniently, we may choose $c$ to be an integer power of $p$ in these cases, so the matrices have entries in $\mathbb{Q}_p$. 
Suppose $2j \ge i > j \ge 0$. Then we shall show the stronger statement that $a_{ij}^{(7/12)}/j = b_{ji}^{(5/12)}/i \in \mathbb{Z}_2$. From \cite{BCslopes} we know that
\[a_{ij}^{(\frac{7}{12})} = 2^{j-i} a_{ij}^{(\frac{1}{2})} =
2^{j-i}\cdot 6ij\left( \frac{(2j)!}{2^j j!}\right)^2 \left( \frac{2^i i!}{(2i)!}\right)^2 \frac{(2i-1)!}{(i+j)!} \frac{(2j+i-1)!}{(3j)!} \binom{j}{i-j}.\]
The first two bracketed terms are clearly in $\mathbb{Z}_2^\times$, so we can safely ignore them. If we put $i = j + t$, what is left is
\[2^{1-t} \cdot 3ij \left(\frac{(2j + 2t - 1)!}{(2j + t)!}\right)\left( \frac{(3j+t-1)!}{(3j)!}\right) \binom{j}{t}.\]
If $t$ is odd, we are safe, as the two factorial terms each simplify to products of $t-1$ consecutive integers, and each product contains $\frac{t-1}{2}$ even integers which cancel all the factors of 2 in the denominator. If $t$ is even, then we are in slightly more trouble. The first product always ends on an odd integer so it has $\frac{t}{2}-1$ even terms, and the second one depends on $j$; if $3j+1$ is even, we get $\frac{t}{2}$ even factors, but if $(3j+1)$ is odd, then we are one short. However, this occurs only if $j$ is even, and consequently $i$ is even; so $a_{ij}/j \in \mathbb{Z}_2$, as claimed.
\end{proof}

\begin{theorem}\label{mainlemma} 
Let $K$ be a field complete with respect to a non-archimedean valuation, with ring of integers $\mathcal{O}_K$ and maximal ideal $\mathfrak{M}_K$. Let $S$ be the space of sequences over $K$ with entries tending to zero. Then if $M$ is any operator on $S$ given by a matrix of the form $ADB$ where $D$ is diagonal with strictly increasing valuations and $A$,$B$ have entries in $\mathcal{O}_K$ congruent to the identity modulo $\mathfrak{M}_K$, then we can find a matrix $C$, also with integral entries congruent to the identity, such that $C^{-1}MC$ is diagonal.
\end{theorem}

\begin{proof} The statement is not affected by conjugating $M$ by any matrix congruent to the identity, so we conjugate by $B^{-1}$, allowing us to assume without loss of generality that $M = AD$.
It is known (see \cite{BCslopes}) that $M$ has the same Newton polygon as $D$. Hence, for every $j$ there is an eigenvector $v_j$ such that $M v_j = \mu_j v_j$ with $\frac{\mu_j}{D_{jj}} \in \mathcal{O}_K^\times$, and $v_j$ is unique up to scalars. We normalise $v_j$ so it is integral with norm $1$.

Suppose $D v_j = \eta_j w_j$, where $w_j$ has norm 1 and $\eta_j \in K$. Then since $A = \mathrm{Id} \bmod \mathfrak{M}_K$, $\mu_j v_j = AD v_j = \eta_j A w_j$. Comparing norms, we see that $\varepsilon_j = \eta_j^{-1} \mu_j \in \mathcal{O}_K^\times$, and reducing mod $\mathfrak{M}_K$ we have $\overline{\varepsilon_j}\  \overline{v_j} = \overline A\ \overline{w_j}$. But $\overline{A}$ is the identity, and consequently $\overline{\varepsilon_j}\ \overline v_j = \overline w_j$. This is impossible unless $\overline{v_j}$ has all its components zero outside the $j$th.

Now if $C$ is the matrix whose $j$th column is $v_j$, then we evidently have $MC = CE$ where $E$ is the diagonal matrix with $E_{ii} = \mu_i$, and since $C$ is congruent to the identity, it is necessarily invertible (since the series $(1 + T)^{-1} = 1 - T + T^2 + \dots$ converges whenever $|T| < 1$).
\end{proof}

\begin{corollary}[Spectral expansion theorem] For any $r \in \left(\frac{5}{12}, \frac{7}{12}\right)$, the finite slope eigenfunctions form an orthonormal basis of the space $\mathcal{S}_0(r)$; that is, for all $h \in \mathcal{S}_0(r)$, the sum \[\sum_{i=1}^\infty \pi_i(h)\]
converges to $h$, and $\|h\| = \sup_{i} \|\pi_i(h)\|$.
\end{corollary}

Note in particular that this implies that the kernel of $U$ is zero for all $r > \frac{5}{12}$; it is in fact known that the kernel is zero for $r \ge \frac{1}{p + 1}$, by Lemma 6.13 of \cite{BCproper}.

\section{Appendix A: Overconvergent forms at small level}
\label{appendix}

In this appendix, we finish off the proof of Theorem \ref{basis} in order to show that the space we work with really is the same as the space of $r$-overconvergent $p$-adic modular forms, for each $p \in \{2, 3, 5, 7, 13\}$. Since we work only with weight zero forms, the problem of whether or not the sheaf $\omega^{\otimes k}$ descends does not arise, and hence the problem is reduced to identifying in terms of our chosen uniformiser the region of $X_0(p)$ corresponding to the $r$-overconvergent locus. For $p \ge 5$, the Hasse invariant lifts to level 1 via the classical level 1 Eisenstein series $E_{p-1}$, so we can measure overconvergence directly using this form; the argument is given in \cite[Prop 3.5]{smithline-bounding}. However, for $p = 2$ and $p = 3$, the Hasse invariant does not lift to characteristic 0 in level 1, so we need to introduce auxiliary level structure. The case $p = 2$ is covered in \cite[\S 4]{BCslopes}, using a weight 1 $\theta$ series of level 3 as a Hasse lifting, so we are left with the case $p = 3$. Smithline shows that in this case the region where $|f_3| \le 3^{6r}$ coincides with the region where $|E_6| \ge 3^{-3r}$, for all $r < \frac{3}{4}$; so we must compare the valuations of $E_6$ and the Hasse invariant. 

Consider the $2$-stabilised Eisenstein series $E_2' = 2E_2(2z)-E_2(z)$, which is a modular form of weight 2 and level $\Gamma_0(2)$. Since $E_2(z) \equiv E_2(2z) \equiv 1 \bmod 3$, $E_2'$ is a lift of the mod 3 Hasse invariant. Using our parameter $f_2$ on $X_0(2)$, we have the identities
\[\frac{E_2'^6}\Delta = \frac{(1 + 2^6 f_2)^3}{f_2}\]
and
\[\frac{E_6^2}\Delta = \frac{(1 + 2^6 f_2)(1 - 2^9 f_2)^2}{f_2}.\]
The supersingular region corresponds to $|1+2^6 f_2| < 1$; in this region $|f_2| = 1$, so if $|1+2^6 f_2| > 3^{-2}$, then $|1+2^6 f_2| = |1 + 2^6f_2 - 9.2^6f_2| = |1 - 2^9f_2|$. Since supersingular curves have good reduction, $|\Delta| = 1$ also, hence 
\begin{align*}
|E_2'| \ge 3^{-r} &\iff \left|\frac{E_2'^6}{\Delta}\right| \ge 3^{-6r} \\
&\iff \left|\frac{E_6^2}{\Delta}\right| \ge 3^{-6r} \\
&\iff |E_6| \ge 3^{-3r}\\
\end{align*}
for all $r < 1$, and the result follows.

\section{Appendix B: $q$-expansions of small slope 5-adic eigenfunctions}
\label{tables}

The following list gives the first 20 terms of the $q$-expansions of the 20 smallest slope $5$-adic eigenforms, with the coefficients given to a relative precision of $O(5^{10})$. This computation took less than 1 minute on a standard laptop PC.

\begin{multline*} \phi_{1} = q + 8528631q^{2} + 8596652q^{3} + 2788848q^{4} + 5\times610813q^{5} + 6727787q^{6}\\
 + 2747331q^{7} + 5\times3412617q^{8} + 6989312q^{9} + 5\times4155753q^{10} + 538817q^{11}\\
 + 9643146q^{12} + 6371187q^{13} + 5536986q^{14} + 5\times9298076q^{15} + 8198461q^{16}\\
 + 3226656q^{17} + 5179372q^{18} + 5\times9335108q^{19} + 5\times7582174q^{20} +  O(q^{21})
\end{multline*}
\begin{multline*} \phi_{2} = q + 441709q^{2} + 2550713q^{3} + 4301618q^{4} + 5^4\times2356503q^{5} + 2966642q^{6}\\
 + 3223594q^{7} + 5\times9703174q^{8} + 7251077q^{9} + 5^4\times9677377q^{10} + 3828592q^{11}\\
 + 5453634q^{12} + 4410268q^{13} + 3763396q^{14} + 5^4\times1117889q^{15} + 1692896q^{16}\\
 + 2395464q^{17} + 4642468q^{18} + 5\times2705229q^{19} + 5^4\times8143729q^{20} + O(q^{21})
\end{multline*}
\begin{multline*} \phi_{3} = q + 7123391q^{2} + 727387q^{3} + 8909193q^{4} + 5^5\times6386403q^{5} + 6931192q^{6}\\
 + 3140781q^{7} + 5\times2842166q^{8} + 3306102q^{9} + 5^5\times3855698q^{10} + 1486467q^{11}\\
 + 1481191q^{12} + 909182q^{13} + 3295871q^{14} + 5^5\times5659586q^{15} + 2077746q^{16}\\
 + 7148211q^{17} + 2935007q^{18} + 5\times6743039q^{19} + 5^5\times1590279q^{20}
 + O(q^{21})
\end{multline*}
\begin{multline*} \phi_{4} = q + 2764444q^{2} + 5364423q^{3} + 7074448q^{4} + 5^8\times6938782q^{5} + 8303937q^{6}\\
 + 2059419q^{7} + 5\times5835813q^{8} + 6128137q^{9} + 5^8\times9032833q^{10} +9024817q^{11}\\
 + 9297879q^{12} + 3774838q^{13} + 3966786q^{14} + 5^8\times3159036q^{15} + 908886q^{16}\\
 + 1286194q^{17} + 2888953q^{18} + 5\times3751388q^{19} + 5^8\times5567336q^{20}
 + O(q^{21})
\end{multline*}
\begin{multline*} \phi_{5} = q + 5791436q^{2} + 3059457q^{3} + 3403033q^{4} + 5^9\times8921438q^{5} + 6832127q^{6}\\
 + 3955981q^{7} + 5\times3439059q^{8} + 6952557q^{9} + 5^9\times7517468q^{10} +9760342q^{11}\\
 + 7351831q^{12} + 8002297q^{13} + 231841q^{14} + 5^9\times79791q^{15} + 4456166q^{16}\\
 + 7616646q^{17} + 5698727q^{18} + 5\times7110866q^{19} + 5^9\times6515204q^{20}
 + O(q^{21})
\end{multline*}
\begin{multline*} \phi_{6} = q + 6831044q^{2} + 1698148q^{3} + 2950248q^{4} + 5^{10}\times6825297q^{5} + 6519012q^{6}\\
 + 8819044q^{7} + 5\times5659178q^{8} + 8713237q^{9} + 5^{10}\times7635693q^{10} + 4926567q^{11}\\
 + 6568829q^{12} + 5335163q^{13} + 6117561q^{14} + 5^{10}\times3121831q^{15} + 9149661q^{16}\\
 + 3456869q^{17} + 7282553q^{18} + 5\times82178q^{19} + 5^{10}\times464281q^{20} + O(q^{21})
\end{multline*}
\begin{multline*} \phi_{7} = q + 8461691q^{2} + 7744062q^{3} + 4618543q^{4} + 5^{13}\times9616002q^{5} + 8166342q^{6}\\
 + 9150156q^{7} + 5\times7971386q^{8} + 1468177q^{9} + 5^{13}\times860632q^{10} +5105092q^{11}\\
 + 4044791q^{12} + 5464782q^{13} + 1658171q^{14} + 5^{13}\times1617624q^{15} + 6957796q^{16}\\
 + 2187611q^{17} + 8154182q^{18} + 5\times4201019q^{19} + 5^{13}\times4662586q^{20} + O(q^{21})
\end{multline*}
\begin{multline*} \phi_{8} = q + 9458634q^{2} + 1415388q^{3} + 310018q^{4} + 5^{14}\times7929152q^{5} + 341242q^{6}\\
 + 8941094q^{7} + 5\times5522594q^{8} + 6133252q^{9} + 5^{14}\times1385868q^{10} + 1356842q^{11}\\
 + 6694484q^{12} + 1201868q^{13} + 8361846q^{14} + 5^{14}\times4325351q^{15} + 165471q^{16}\\
 + 8543864q^{17} + 8163393q^{18} + 5\times8748199q^{19} + 5^{14}\times6016611q^{20} + O(q^{21})
\end{multline*}
\begin{multline*} \phi_{9} = q + 1036606q^{2} + 8499877q^{3} + 6100798q^{4} + 5^{19}\times9288232q^{5} + 7872462q^{6}\\
 + 6770081q^{7} + 5\times8252407q^{8} + 2114087q^{9} + 5^{19}\times4598717q^{10} + 7406442q^{11}\\
 + 7211221q^{12} + 9554887q^{13} + 6194461q^{14} + 5^{19}\times1422464q^{15} + 9065311q^{16}\\
 + 5385831q^{17} + 659347q^{18} + 5\times9351018q^{19} + 5^{19}\times2209136q^{20}
 + O(q^{21})
\end{multline*}
\begin{multline*} \phi_{10} = q + 8935814q^{2} + 2184043q^{3} + 7194158q^{4} + 5^{20}\times9176128q^{5} + 844127q^{6}\\
 + 1292144q^{7} + 5\times1755091q^{8} + 8018557q^{9} + 5^{20}\times1173192q^{10} + 9267217q^{11}\\
 + 6670794q^{12} + 8784078q^{13} + 1023341q^{14} + 5^{20}\times3438004q^{15} + 9735791q^{16}\\
 + 7839479q^{17} + 9681648q^{18} + 5\times9158266q^{19} + 5^{20}\times2941474q^{20} + O(q^{21})
\end{multline*}
\begin{multline*} \phi_{11} = q + 8097156q^{2} + 5482427q^{3} + 4624273q^{4} + 5^{21}\times3090372q^{5} + 6130737q^{6}\\
 + 9435206q^{7} + 5\times3663802q^{8} + 7112412q^{9} + 5^{21}\times1525782q^{10} + 9588067q^{11}\\
 + 2822446q^{12} + 9371737q^{13} + 4796011q^{14} + 5^{21}\times4517844q^{15} + 9306236q^{16}\\
 + 2578856q^{17} + 6765897q^{18} + 5\times5575723q^{19} + 5^{21}\times1143306q^{20} + O(q^{21})
\end{multline*}
\begin{multline*} \phi_{12} = q + 9675784q^{2} + 6753913q^{3} + 116218q^{4} + 5^{24}\times8946888q^{5} + 8811542q^{6}\\
 + 8069219q^{7} + 5\times6507279q^{8} + 2269902q^{9} + 5^{24}\times1463317q^{10} + 3569092q^{11}\\
 + 4386034q^{12} + 8715668q^{13} + 4467696q^{14} + 5^{24}\times9032119q^{15} + 3147446q^{16}\\
 + 1255689q^{17} + 5281293q^{18} + 5\times2446659q^{19} + 5^{24}\times4273334q^{20} + O(q^{21})
\end{multline*}
\begin{multline*} \phi_{13} = q + 852841q^{2} + 6464712q^{3} + 8669718q^{4} + 5^{25}\times9222513q^{5} + 2306167q^{6}\\
 + 7752656q^{7} + 5\times1741296q^{8} + 4498152q^{9} + 5^{25}\times5178183q^{10} + 8249092q^{11}\\
 + 3428716q^{12} + 4365957q^{13} + 5551946q^{14} + 5^{25}\times1789381q^{15} + 4399821q^{16}\\
 + 7853311q^{17} + 7277957q^{18} + 5\times2773209q^{19} + 5^{25}\times4930084q^{20} + O(q^{21})
\end{multline*}
\begin{multline*} \phi_{14} = q + 3696344q^{2} + 5088573q^{3} + 4864773q^{4} + 5^{28}\times7513547q^{5} + 4948987q^{6}\\
 + 9082919q^{7} + 5\times5387723q^{8} + 4212787q^{9} + 5^{28}\times7887793q^{10} + 5486817q^{11}\\
 + 6445179q^{12} + 264638q^{13} + 9163761q^{14} + 5^{28}\times9742181q^{15} + 5608361q^{16}\\
 + 3782269q^{17} + 5653853q^{18} + 5\times2678998q^{19} + 5^{28}\times1361081q^{20} + O(q^{21})
\end{multline*}
\begin{multline*} \phi_{15} = q + 5997936q^{2} + 2852832q^{3} + 6767908q^{4} + 5^{29}\times3494278q^{5} + 239127q^{6}\\
 + 8242231q^{7} + 5\times5754659q^{8} + 1331682q^{9} + 5^{29}\times8544583q^{10} + 7742217q^{11}\\
 + 9202956q^{12} + 5295922q^{13} + 6847716q^{14} + 5^{29}\times4110921q^{15} + 7441416q^{16}\\
 + 6452396q^{17} + 9423977q^{18} + 5\times2768516q^{19} + 5^{29}\times6717924q^{20} + O(q^{21})
\end{multline*}
\begin{multline*} \phi_{16} = q + 7855519q^{2} + 4239748q^{3} + 3954673q^{4} + 5^{30}\times8731987q^{5} + 1047337q^{6}\\
 + 7593044q^{7} + 5\times6656568q^{8} + 7374337q^{9} + 5^{30}\times2676878q^{10} + 5407692q^{11}\\
 + 536154q^{12} + 2961238q^{13} + 6487961q^{14} + 5^{30}\times8403651q^{15} + 524436q^{16}\\
 + 8063044q^{17} + 6134653q^{18} + 5\times7095743q^{19} + 5^{30}\times4787751q^{20} + O(q^{21})
\end{multline*}
\begin{multline*} \phi_{17} = q + 3058366q^{2} + 808487q^{3} + 3957143q^{4} + 5^{35}\times4332043q^{5} + 2667867q^{6}\\
 + 2677656q^{7} + 5\times9265831q^{8} + 2140627q^{9} + 5^{35}\times8177988q^{10} + 8770592q^{11}\\
 + 1797641q^{12} + 6220257q^{13} + 4023221q^{14} + 5^{35}\times7870816q^{15} + 1693096q^{16}\\
 + 9074636q^{17} + 4429232q^{18} + 5\times7074024q^{19} + 5^{35}\times3867524q^{20} + O(q^{21})
\end{multline*}
\begin{multline*} \phi_{18} = q + 4792184q^{2} + 9735438q^{3} + 3075793q^{4} + 5^{36}\times9618893q^{5} + 6310342q^{6}\\
 + 6556094q^{7} + 5\times1549289q^{8} + 6307052q^{9} + 5^{36}\times7085437q^{10} + 8972592q^{11}\\
 + 2599209q^{12} + 5715468q^{13} + 956796q^{14} + 5^{36}\times5570759q^{15} + 552671q^{16}\\
 + 2538389q^{17} + 2900318q^{18} + 5\times6364319q^{19} + 5^{36}\times1100899q^{20} + O(q^{21})
\end{multline*}
\begin{multline*} \phi_{19} = q + 3408581q^{2} + 217102q^{3} + 1581998q^{4} + 5^{39}\times2535503q^{5} + 9752262q^{6}\\
 + 1937831q^{7} + 5\times7503797q^{8} + 7627362q^{9} + 5^{39}\times6413743q^{10} + 6787817q^{11}\\
 + 7664171q^{12} + 3969712q^{13} + 21561q^{14} + 5^{39}\times3787931q^{15} + 4478661q^{16}\\
 + 4153256q^{17} + 4630822q^{18} + 5\times6899078q^{19} + 5^{39}\times8331244q^{20} + O(q^{21})
\end{multline*}
\begin{multline*} \phi_{20} = q + 7376064q^{2} + 5111168q^{3} + 6655533q^{4} + 5^{40}\times1816457q^{5} + 9314002q^{6}\\
 + 8378394q^{7} + 5\times6422316q^{8} + 8376307q^{9} + 5^{40}\times2303998q^{10} + 9013467q^{11}\\
 + 8230044q^{12} + 8742078q^{13} + 48716q^{14} + 5^{40}\times7907401q^{15} + 463666q^{16}\\
 + 6617104q^{17} + 6593773q^{18} + 5\times2535366q^{19} + 5^{40}\times7084706q^{20} + O(q^{21})
\end{multline*}

\section*{Acknowledgements}

I would like to thank Barry Mazur for initially bringing this problem to my attention, and for many helpful conversations while I was working on it; Kevin Buzzard and Frank Calegari, for assistance with the proof of Theorem \ref{mainlemma}; and finally the anonymous referee, whose suggestions improved the exposition substantially.

\bibliographystyle{../halpha}
\bibliography{../library}

\end{document}